\documentclass[12pt]{article}
\usepackage{epsfig}
\usepackage{amsmath}

\newcommand{\nit}{\noindent}

\newcommand{\be}{\begin{equation}}

\newcommand{\ee}{\end{equation}}

\newcommand{\ba}{\begin{eqnarray}}

\newcommand{\ea}{\end{eqnarray}}

\begin{document}
\title{\hspace{.2 in} \\ 
{\Large \bf An Invertible Discrete Auditory Transform}}

\author{Jack Xin \thanks{Corresponding author, Department of Mathematics and ICES,
University of Texas at Austin,
Austin, TX 78712, USA; email:jxin@math.utexas.edu.} \hspace{.1 in} and 
\hspace{.1 in} Yingyong Qi \thanks{
Qualcomm Inc,
5775 Morehouse Drive,
San Diego, CA 92121, USA.}
}

\date{}

\maketitle

\vspace{0.25 in}
\thispagestyle{empty}



\begin{abstract}
A discrete auditory transform (DAT) from sound signal
to spectrum is presented and shown to be invertible in closed form.
The transform preserves energy, and its spectrum is smoother than that of  
the discrete Fourier transform (DFT) consistent with human audition. 
DAT and DFT are compared in signal denoising tests with spectral thresholding method. The signals 
are noisy speech segments. It is found that DAT can gain 3 to 5 decibel (dB) in signal to noise 
ratio (SNR) over DFT except when the noise level is relatively low. 
\end{abstract}

\vspace{.2 in}

\hspace{.2 in} {\bf Keywords: Invertible Discrete
Auditory Transform.}

\newpage
\setcounter{page}{1}

\section{Introduction}
\setcounter{equation}{0}
Audible acoustic signal processing often consists of frame by frame discrete 
Fourier transform (DFT) of input signal followed by spreading in Fourier 
spectrum using the critical band filter, see M. Schroeder et al \cite{S}. 
These two steps mimic the responses of human audition to the input signal, 
facilitating the computation
of excitation pattern over critical bands and psychoacoustics-based processing
\cite{S}. However, the spreading operation in step two is not invertible.  

In this work, an invertible discrete auditory transform (DAT) is formulated to combine the two 
steps into one. DAT incorporates the spectral 
spreading functions of Schroeder et al \cite{S}, which leads to 
smoother spectrum than that of the discrete Fourier transform (DFT).
As a result, DAT has better localization properties in 
the time domain. DAT bears some resemblance to the wavelet transform \cite{ID} 
in that a function of one variable (time) is transformed into 
a function of two variables (time and frequency). It is this 
redundancy that makes the inversion possible and explicit.

The paper is organized as follows. In section 2, a general 
form of DAT is introduced and its inversion 
established. In section 3, a specific DAT is given based on the known auditory spectral energy 
spreading functions of Schroeder et al \cite{S}.
DAT spectrum is defined and compared with DFT spectrum. Time localization 
property of DAT basis functions is illustrated as well.
In section 4, noisy signal reconstruction from thresholded spectrum is carried out and its 
signal to noise ratio is computed by using the reconstructed signal and the clean (noise free) signal. 
The signal is a segment of male or female speech, and 
can be either voiced (e.g. vowels) or unvoiced (consonants such as s and f). 
DAT is found to gain by 3 to 5 decibel (dB) in signal to noise ratio (SNR) over DFT in such a denoising task.
Concluding remarks are made in section 5.

\section{Discrete Auditory Transform}
\setcounter{equation}{0}
Let $s=(s_0,\cdots,,s_{N-1})$ be a discrete signal, and 
$\hat{s}$ its discrete Fourier transform (DFT) \cite{Br}:
\be
 \hat{s}_{k} = \sum_{n=0}^{N-1}\, s_{n}\, e^{-i (2\pi n k/N )}. \label{dt1}
\ee
The DFT inversion formula is:
\be
s_{n} = {1\over N} \, \sum_{k=0}^{N-1}\, \hat{s}_{k}\, e^{i (2\pi n k/N )}. 
\label{dt2}
\ee
For two signals $s$ and $t$ of length $N$, the Plancherel-Parseval equality is:
\be
\sum_{n=0}^{N-1}\, s_{n}t_{n}^{*} = {1\over N}\sum_{k=0}^{N-1}\, 
\hat{s}_{k}\hat{t}_{k}^{*}, \label{dt4}
\ee
implying the energy identity: 
\be \sum_{n=0}^{N-1}\, |s_n|^2 = {1\over N} \sum_{k=0}^{N-1}\, |\hat{s}_{k}|^{2}. \label{df5}
\ee

Define the discrete auditory transform (DAT):
\be
S_{j,m} \equiv \sum_{l=0}^{N-1}\, s_l \, K_{j-l,m}, \label{dt5}
\ee
where the double indexed discrete kernel function is given by:
\be
K_{l,m} = \sum_{n=0}^{N-1}\, X_{m,n} \, 
e^{i (2\pi l n/N)}; \label{dt6}
\ee
where the matrix $X_{m,n}$ has square sum equal to one in $m$:
\be
\sum_{m=0}^{M-1}\, |X_{m,n}|^{2}  =1, \;\; \forall n. 
\label{nm}
\ee

\subsection{Energy Identity}
Let us show the energy conservation property of the transform.
Upon substituting (\ref{dt6}) into (\ref{dt5}), the transform can be written as:
\be
S_{j,m}
 =  \sum_{n=0}^{N-1}\, \hat{s}_{n}\, 
X_{m,n} \, e^{i(2\pi n j/N)}, \label{dt8}
\ee
which is similar to the representation of time domain solutions of cochlear models 
as a sum of time harmonic solutions \cite{xq04,xqd}.

It follows from (\ref{dt8}), (\ref{df5}), and (\ref{dt2}) that:
\[
{1\over N} \sum_{j=0}^{N-1}\, |S_{j,m}|^{2} = \sum_{n=0}^{N-1}\, 
|\hat{s}_{n}|^2 \, |X_{m,n}|^2, \]
implying:
\be
{1\over N^2}\sum_{m=0}^{M-1}\,\sum_{j=0}^{N-1} |S_{j,m}|^2 
= {1\over N} \sum_{n=0}^{N-1}\, |\hat{s}_{n}|^2 = 
\sum_{j=0}^{N-1}\, |s_j|^2. \label{dt9}
\ee
Polarizing with (\ref{dt9}), one finds the analogous Plancherel-Parseval
identity of DAT for two signals $s$ and $t$:
\be
{1\over N^2}\sum_{m=0}^{M-1}\, \sum_{j=0}^{N-1}\, S_{j,m}T_{j,m}^{*} 
= \sum_{j=0}^{N-1}\, s_j t_{j}^{*}. \label{dt10}
\ee

\subsection{Inversion}
The explicit inversion formula is:
\be
s_{j} = {1\over N^2}\sum_{m=0}^{M-1}\,
 \sum_{l=0}^{N-1}\, S_{l,m} \, \sum_{n=0}^{N-1}\, 
X^{*}_{m,n}\, e^{i(2\pi (j -l)n/N)}. \label{dt11}
\ee

\nit {\it Proof:} Consider the sum in $l$. In view of (\ref{dt8}), we see that
\[
{1\over N}\, \sum_{l=0}^{N-1}\, S_{l,m}\, e^{i(2\pi (j -l)n/N)}
= e^{2\pi i j n/N}\, \hat{s}_{n}\, X_{m,n}. \]
So the right hand side of (\ref{dt11}) is equal to:
\[
{1\over N}\, \sum_{m=0}^{M-1}\, \sum_{n=0}^{N-1}\, 
|X_{m,n}|^2\, e^{2\pi i j n/N}\, \hat{s}_{n}, \]
which equals upon summing over $m$ and using (\ref{nm}):
\[ {1\over N}\, \sum_{n=0}^{N-1}\, 
 e^{2\pi i j n/N}\, \hat{s}_{n} =s_{j}. \]
 
\section{Transform Kernel and Spectrum}
\setcounter{equation}{0}
The role of the transform kernel $X_{mn}$ is to 
spread the DFT vector $\hat{s}_{n}$. Here our knowledge of human audition 
will be utilized. Let us adopt the real nonnegative energy spreading function of Schroeder 
et al \cite{S}, denoted by $S_{m,n}=S(b(f_m),b(f_n))$, where $f_m$ is the frequency 
to spread from, $f_n$ is the frequency to spread to, and $b$ is the standard mapping from 
Hertz (Hz) to Bark scale \cite{Hart}. The functional form of $S(\cdot,\cdot)$ is given in 
Schroeder et al \cite{S}. 
    
The DFT of a real vector $s$ satisfies the symmetry property 
$\hat{s}_{k} =\hat{s}^{*}_{N-k}$, $k=1,2\cdots, N-1$. 
It is natural for the spreading kernel $X_{m,n}$ to respect this 
symmetry. Suppose the discrete signal $s$ 
has sampling frequency $F_s$ (Hz). The DFT component $\hat{s}_{n}$ 
($0 \leq n \leq N/2$, $N$ a power of 2) corresponds to frequency:
\be
f_n = F_s \cdot n/N, \;\;\;\; n \leq N/2.   \label{dt16}
\ee
Let $U_{m,n}=U(f_m,f_n) = S^{1/2}(b(f_m),b(f_n))$, $0 \leq m \leq M-1$, $M= N/2$. The square 
root is to convert spreading from energy to amplitude scale. Then normalize $U_{m,n}$ to 
define $X_{m,n}$ as follows:

\[X_{m,0} ={U(f_m,f_1) \over m_f (f_1)}, \]

\[X_{m,n}={U(f_m,f_n) \over m_f (f_n)}, \;\;\; \; 1\leq n \leq N/2 -1,\]

\[X_{m,n}={U(f_m,f_{N-n}) \over m_f (f_{N-n})},\;\;\;\; 
N/2 \leq n \leq N-1,\]
where the $m_f$ function is:
\be
m_f(f) =\left( \sum_{m=0}^{M-1} \, |U(f_m,f)|^{2}\right )^{1/2}. \label{dt7}
\ee
We see that $X_{m,n}$ is symmetric in $n$
with respect to $N/2$, and periodic in $n$. The normalization property
(square sum equal to one) with respect to $m$ holds. 
See Figure 1 for a plot of $X_{m,n}$ in $n$, at $m=5:10:55$, where 
$N=128$, $Fs =16000$ Hz. 

\subsection{DAT Spectrum}
In view of formula (\ref{dt8}), we extract the amplitude (intensity) of the 
$S_{j,m}$ ($j\in [0,N-1]$) for each $m$ and define the DAT spectrum:
\be
{\rm spec}(m) \equiv \left ( \sum_{n=0}^{N-1}\, \left | \hat{s}_{n}\, X_{m,n} \right |^{2}
\right )^{1/2}. \label{dt18}
\ee
As the input sound signal is divided into 
frames of length $N$, DAT spectrum can vary from frame to frame in time. 

\subsection{Transform Properties}
Let us illustrate the DAT properties by considering a 500 Hz square wave (top panel of Figure 2).
Middle and bottom panels of Figure 2 show DFT and DAT spectra of 
the 500 Hz square wave in the first frame, for $N=128$, and $M=N/2=64$.  
The later frames are similar. Compared with DFT spectrum, 
DAT spectrum is smoother, especially towards higher frequencies. 

In the time domain representation (\ref{dt11}), $N$ times the inverse DFT
of $X^{*}_{m,n}$ in $n$ for each $m$ plays the role of 
basis functions. Figure 3 shows such functions for $m=20$ and $m=40$, 
their time domain localization property reflects the smoother DAT spectrum.

\section{DAT and DFT in Signal Denoising}
DAT and DFT were used to denoise speech signals. Both were numerically implemented with FFT. 
A simple thresholding method in the transformed domain was applied to improve the 
signal-to-noise ratio (SNR) of noisy speech. The underlying assumption of the method 
is that low level components in the transformed domain are more likely to be noise 
than signal plus noise. Thresholding, therefore, could improve the overall SNR of 
the signal. It is a simple denoise method. More advanced methods exist for 
noise reduction and will be studied in the future. 
The simple thresholding method serves as a tool here to 
reveal the difference between DAT and DFT in signal processing. 

Voiced and unvoiced speech segments were selected from a male and a female speaker
respectively. Each segment has 512 data points. Noisy speech was created by adding 
Gaussian noise to the selected segments. The level of noise was set to produce the 
SNR ranging from -12 dB to +12 dB with a 3 dB step size. DAT and DFT were applied 
to the noisy speech signals. The magnitude of transformed components were then 
compared to a threshold. All components with magnitude smaller than the threshold 
were ignored for the reconstruction of the signal. Here, the threshold was computed 
as the average of the DFT magnitude spectrum. Signal was reconstructed directly by 
the inverse DAT and DFT, respectively. The SNRs of the reconstructed signal was computed and shown 
in Figure 4. Samples of original signal, its noise added signal  (SNR = 0 dB), and 
the signal denoised by thresholding   
were shown for a voiced (Figure 5) and unvoiced (Figure 6) speech segment. 

These results indicate that DAT thresholding has about 3 to 5 dB SNR gain 
over DFT thresholding for voiced speech signals. The improvement is larger 
when the noise level is relatively high. For unvoiced (noise-like) speech signal, 
DAT thresholding also has a SNR gain when the noise level is relatively high. 
The DFT thresholding, however, has higher SNRs when noise level is relatively low.
The DAT thresholding appears to have difficulty in discriminating between 
the original and added noise for unvoiced speech segments when the 
noise added is relatively low. This should not be as handicapping, 
in part because it may not be as necessary to denoise low-level 
noise when the speech itself is noise-like (see Figure 6). Denoise is mostly 
needed for voiced signals with high-level noise. The  potential 
advantage of DAT thresholding is well demonstrated by the 3 to 5 dB 
improvement of SNR when noise level is relatively high. 
The noise-reduction advantage is likely a result of the spectral 
spreading (weighted local spectral averaging) operation of DAT.

\section{Concluding Remarks}
Discrete auditory transforms (DAT) are introduced and shown to be invertible and 
energy preserving. DAT spectra are smoother than DFT's, and DAT basis functions are 
more localized than DFT's in the time domain. 
In signal denoising with spectral thresholding method, it is observed that 
DAT increased SNR by 3 to 5 dB over DFT.
Further study of DAT will be worthwhile in more complicated signal 
processing tasks.

\section{Acknowledgements}
This work was supported in part by NSF grant ITR-0219004,  
a Fellowship from the John Simon Guggenheim Memorial Foundation,  
and the Faculty Research Assignment Award at the University of Texas at Austin.
We thank Prof. G. Papanicolaou and Prof. S-T Yau for their interest and helpful comments.  

\vspace{.2 in}

\bibliographystyle{plain}

\begin{thebibliography}{99}
\vspace{.1 in}

\bibitem{Br}P. Br\'emaud,
``Mathematical Principles of Signal Processing: Fourier and Wavelet
Analysis'', Springer-Verlag, 2002.

\bibitem{ID}I. Debauchies,
``Ten Lectures on Wavelets'', CMS-NSF Regional Conference in
Applied Mathematics, SIAM, Philadelphia, 1992.
                                                    
\bibitem{Hart}W. M. Hartmann,
``Signals, Sound, and Sensation'', Springer, 2000, pp 251-254.

\bibitem{S}M. R. Schroeder, B. S. Atal and J. L. Hall,
{\em Optimizing digital speech
coders by exploiting properties of the human ear},
Journal Acoust. Soc. America,  66(6), pp 1647-1652 (1979).

\bibitem{xq04}J. Xin and Y. Qi,
{\em Global well-posedness and multi-tone solutions
of a class of nonlinear nonlocal cochlear models in hearing}, 
Nonlinearity 17(2004), pp 711-728.

\bibitem{xqd}J. Xin, Y. Qi, and L. Deng,
{\em Time domain computation of a nonlinear nonlocal cochlear model 
with applications to multitone interaction in hearing},
Comm. Math. Sci., Vol. 1, No. 2, 2003, pp. 211--227.

\end{thebibliography}










\newpage
\begin{figure}[p]
\centerline{\includegraphics[width=350pt,height=340pt]{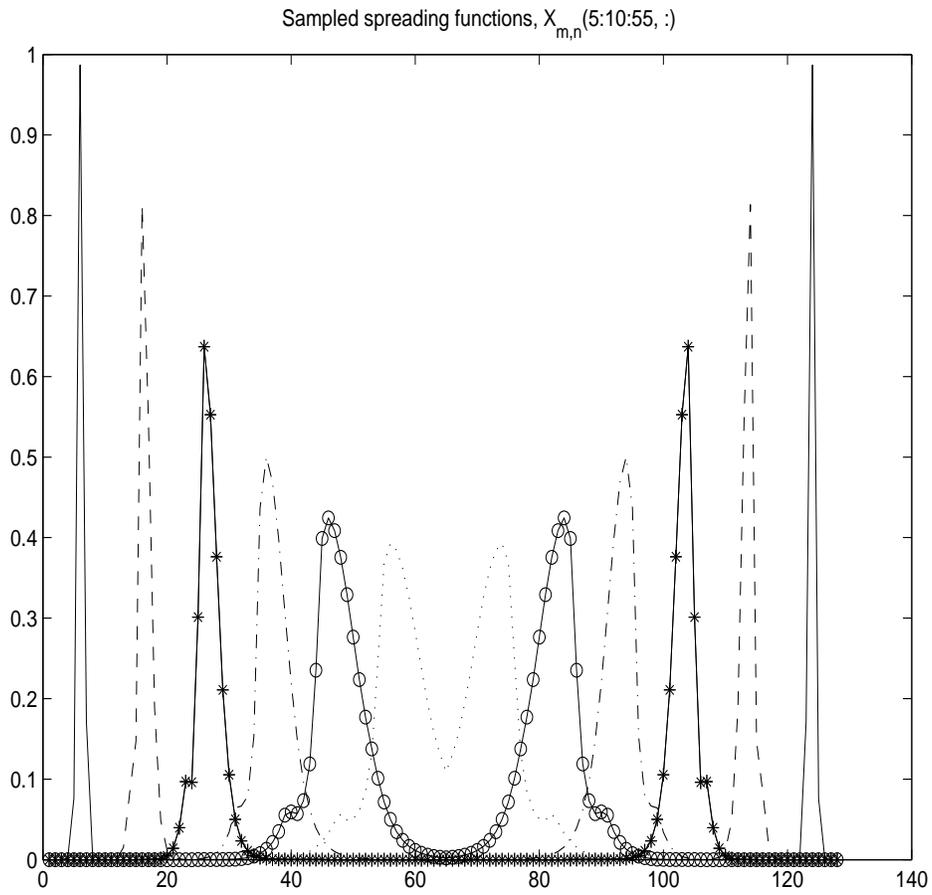}}
\vspace{.25 in}
\caption{The transform kernel $X_{m,n}$ as a function of integer $n =$ 1 to 128. 
The $m$ variable equals 5 (solid), 15 (dash), 25 (line-star), 35 (dashdot), 45 (line-circle), 55 (dot).}
\end{figure}

\newpage
\begin{figure}
\centerline{\includegraphics[width=310pt,height=100pt]{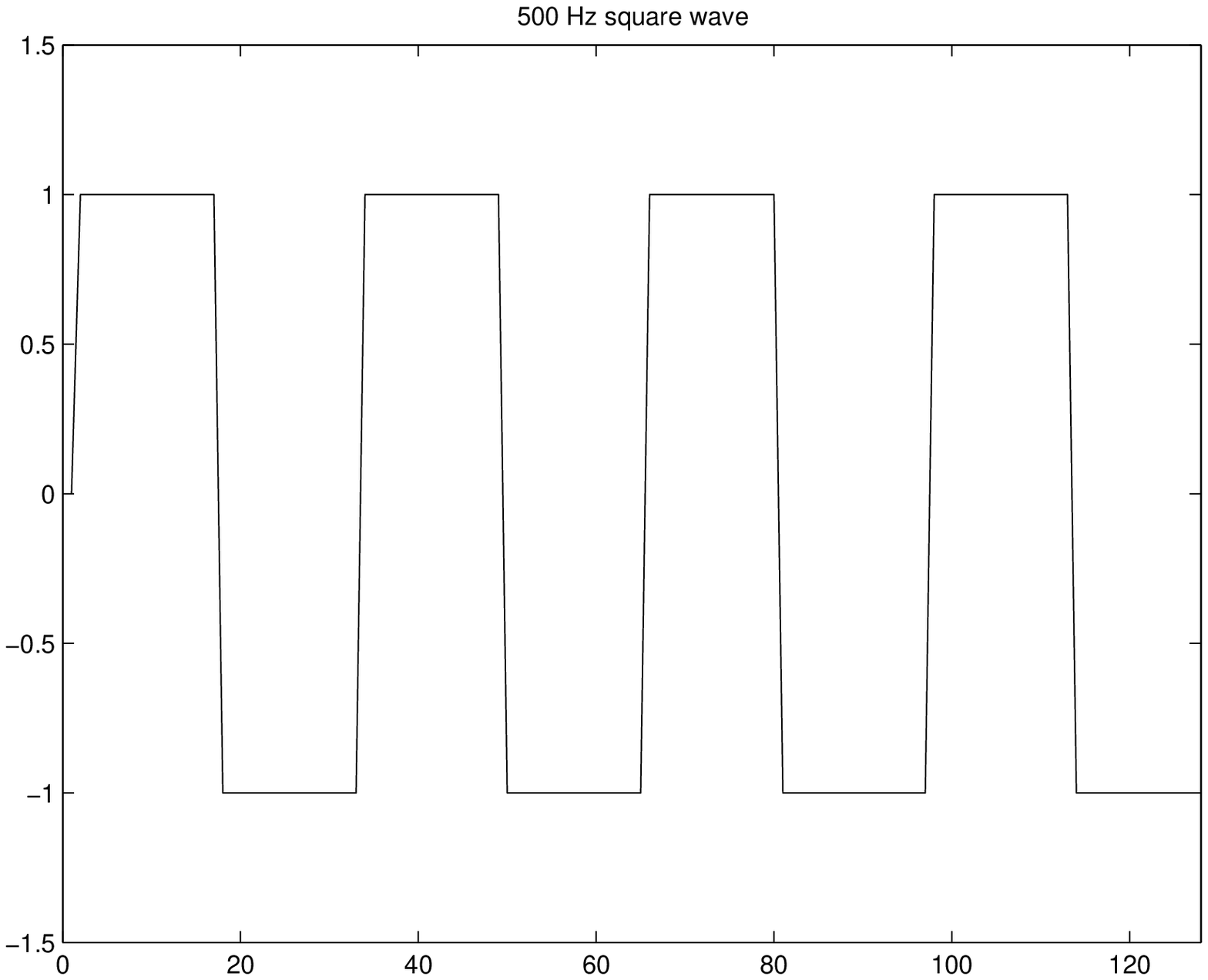}}
\vspace{.1 in}
\centerline{\includegraphics[width=350pt,height=240pt]{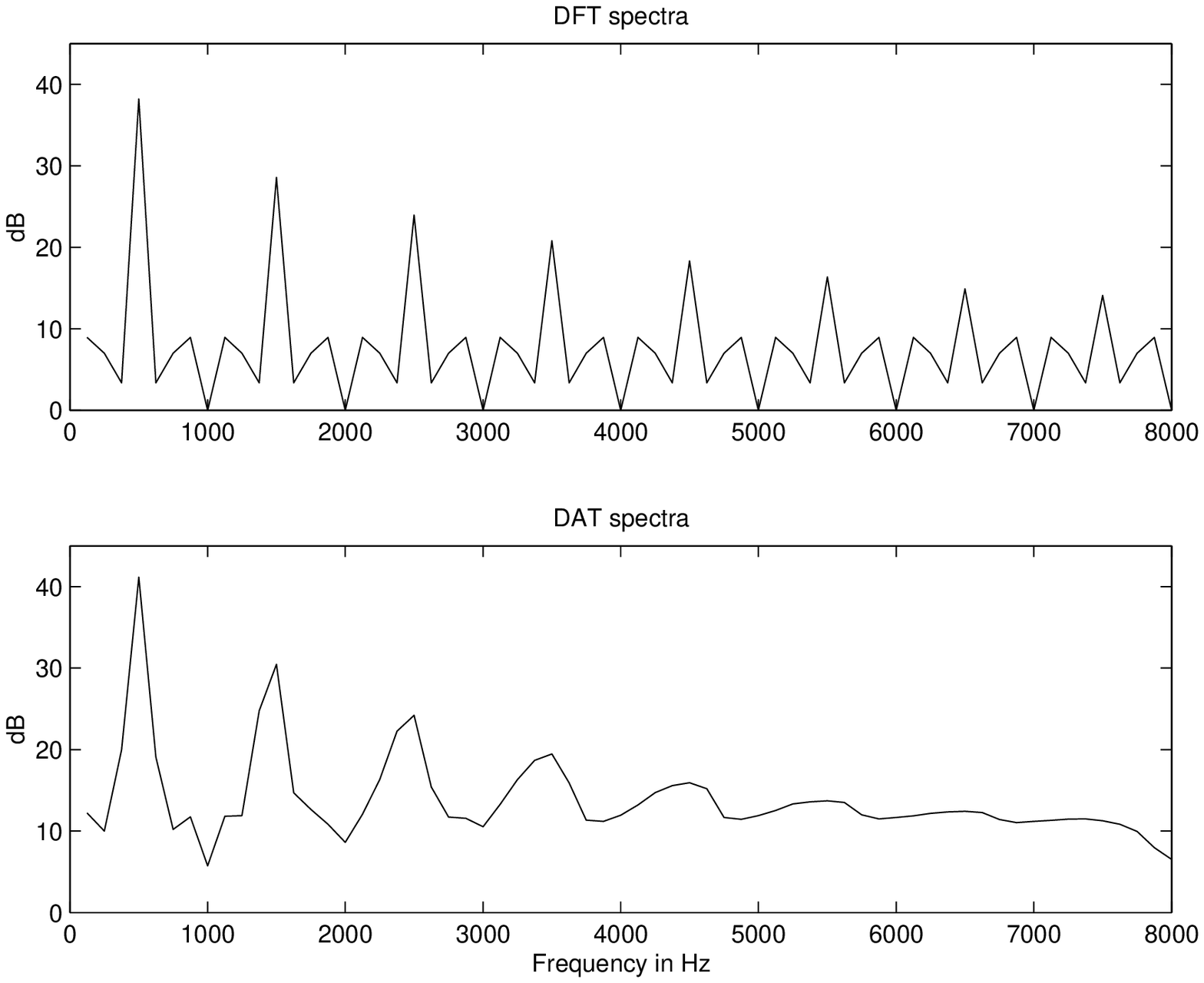}}
\vspace{.25 in}
\caption{A 500 Hz square wave (top), its DFT spectrum (middle) and DAT spectrum (bottom).}
\end{figure}

\newpage
\begin{figure}
\centerline{\includegraphics[width=350pt,height=340pt]{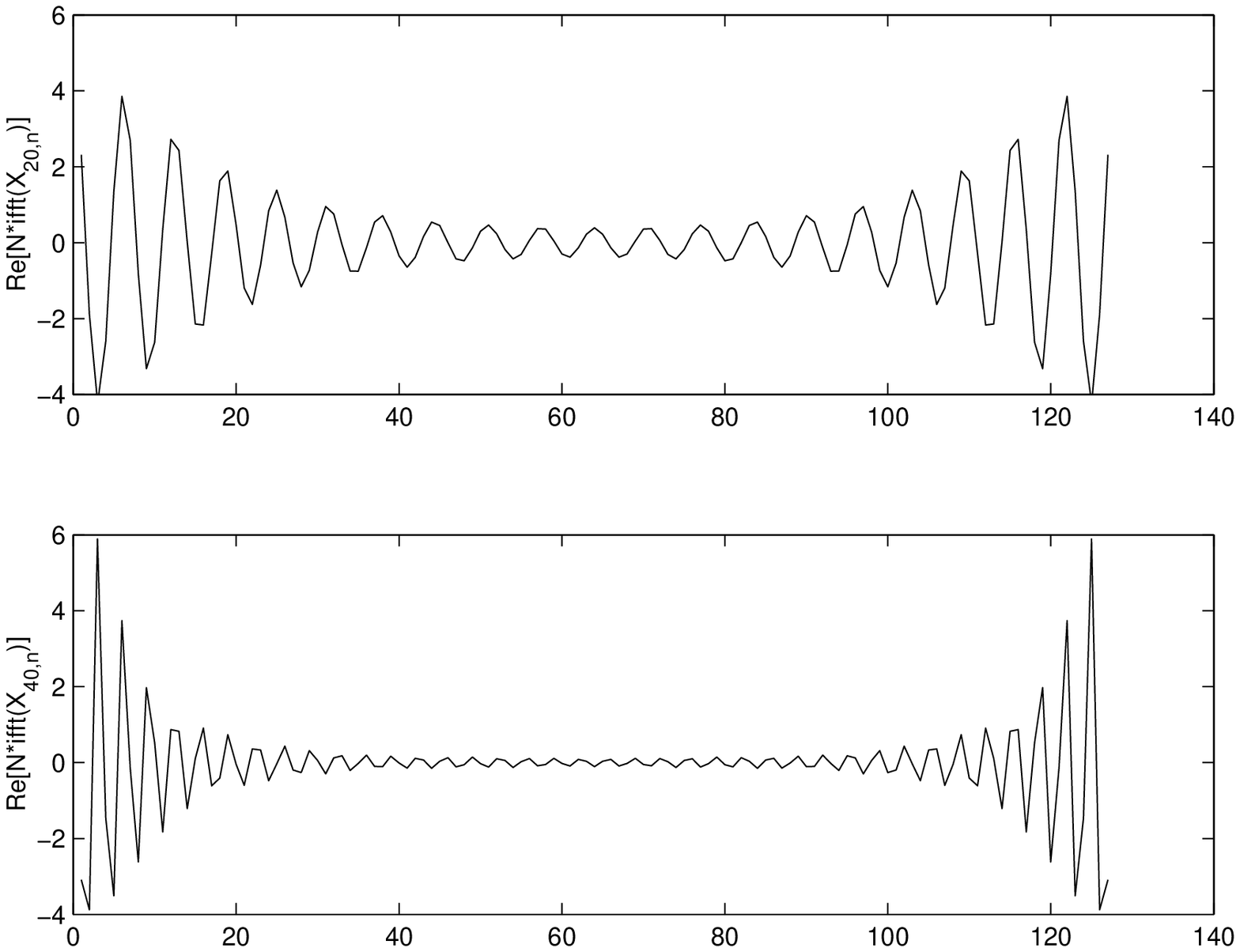}}
\vspace{.25 in}
\caption{An illustration of localization property of DAT basis functions in the time domain.}
\end{figure}


\newpage
\begin{figure}
\centerline{\includegraphics[width=480pt,height=480pt]{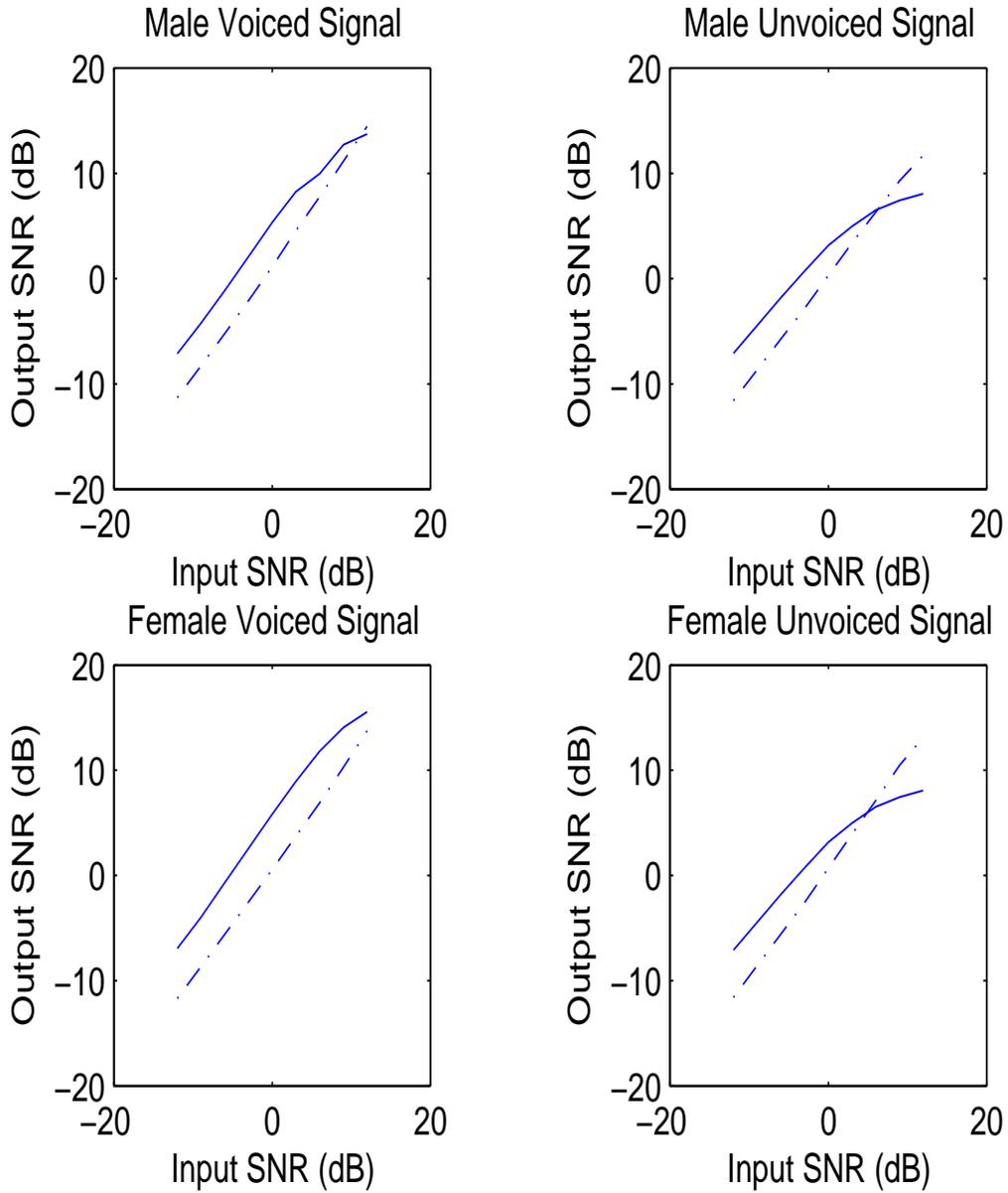}}
\vspace{.25 in}
\caption{Comparison of DAT (solid) and DFT (dashdot) denoising by spectral thresholding for 
male/female, voiced/unvoiced speech segments.}
\end{figure}

\newpage
\begin{figure}
\centerline{\includegraphics[width=350pt,height=350pt]{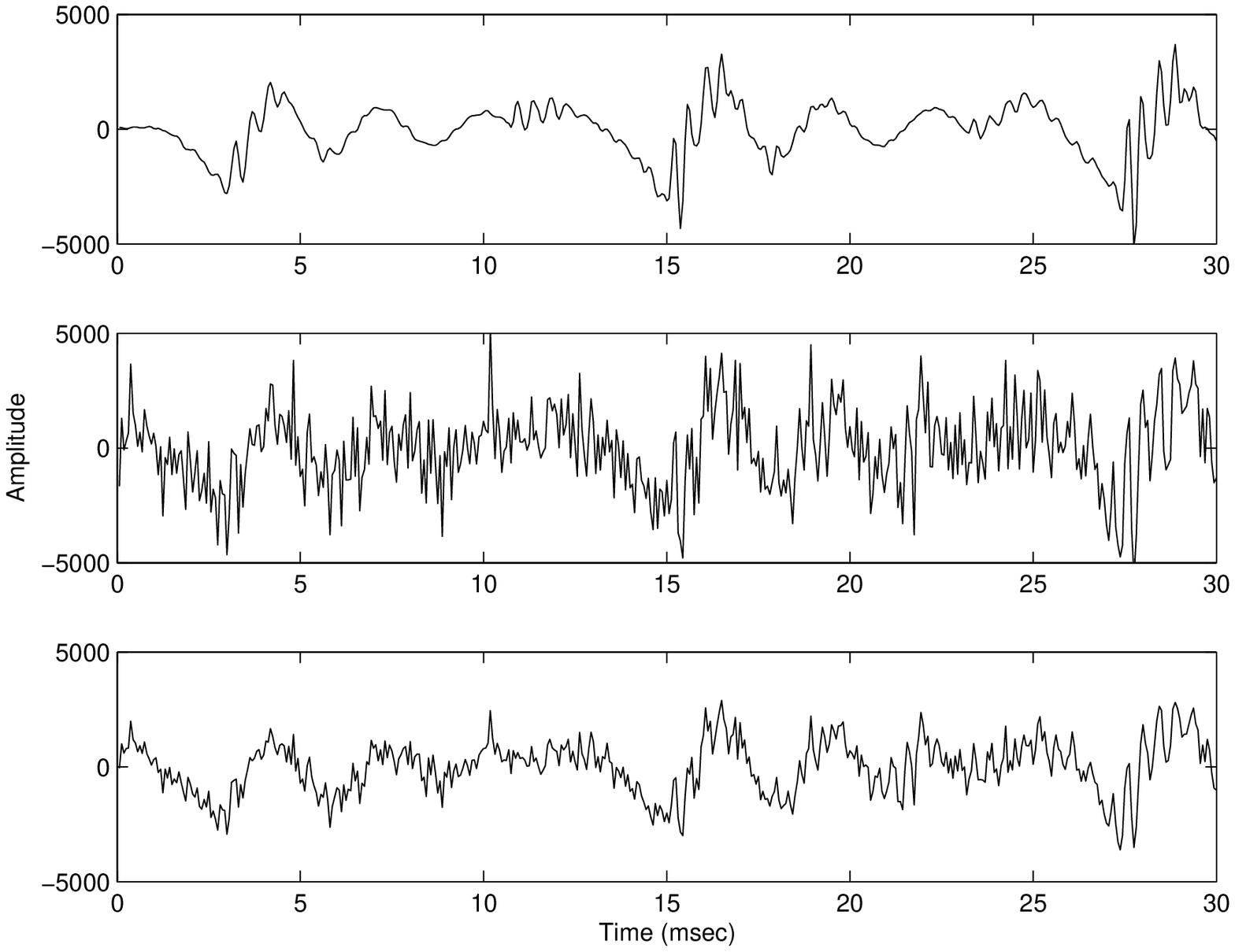}}
\vspace{.25 in}
\caption{(Top down) voiced speech signal, noisy signal (SNR = 0 dB), denoised signal.}
\end{figure}

\newpage
\begin{figure}
\centerline{\includegraphics[width=350pt,height=350pt]{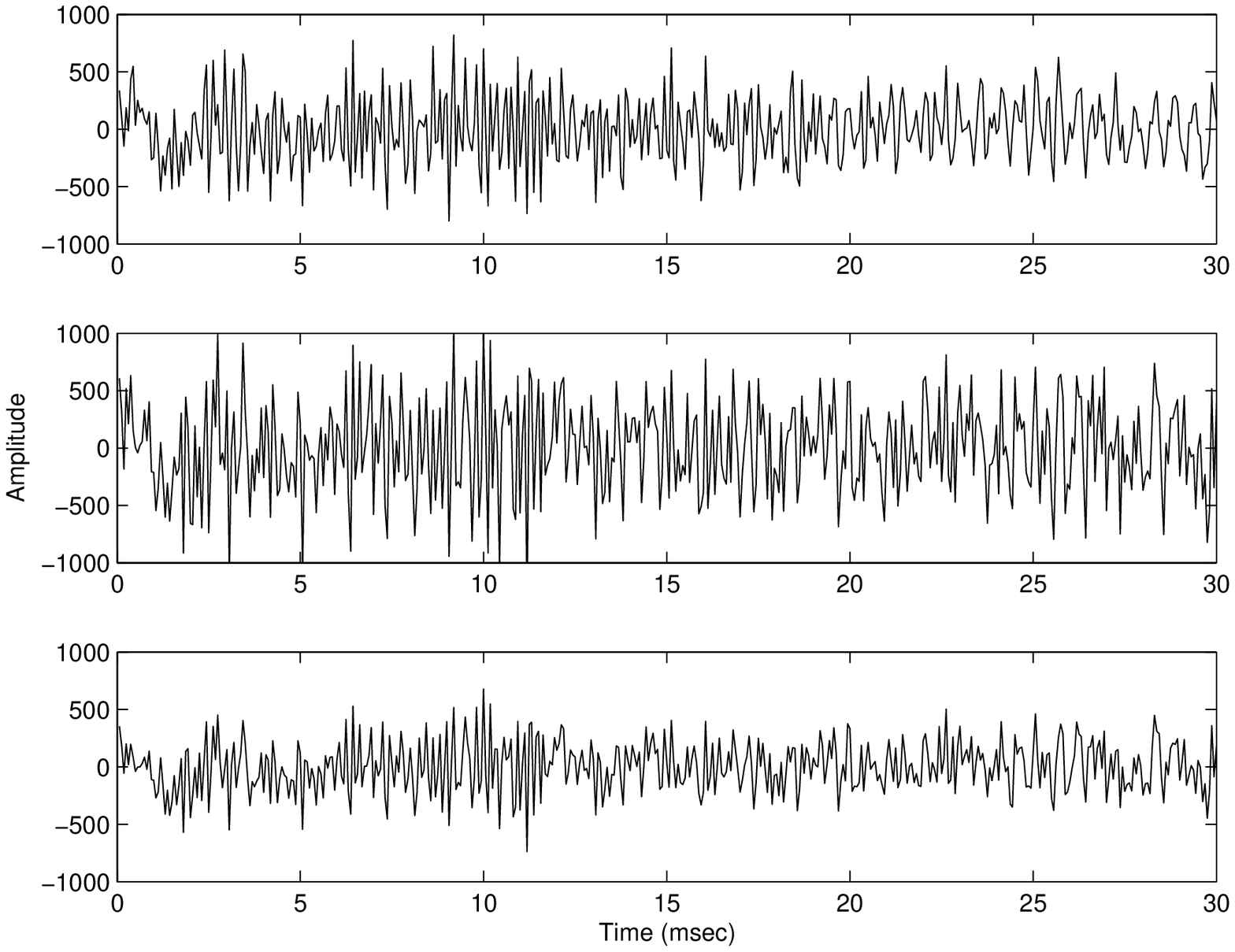}}
\vspace{.25 in}
\caption{(Top down) unvoiced speech signal, noisy signal (SNR = 0 dB), denoised signal.}
\end{figure}

\end{document}